\documentclass[12pt]{article}
\usepackage{amssymb}
\begin{document}

\newcommand{\nc}[2]{\newcommand{#1}{#2}}
\newcommand{\rnc}[2]{\renewcommand{#1}{#2}}
\nc{\bpf}{{\it Proof.~~}}
\nc{\bpr}{\begin{prop}}
\nc{\bth}{\begin{theorem}}
\nc{\ble}{\begin{lem}}
\nc{\bco}{\begin{corollary}}
\nc{\bre}{\begin{remark}}
\nc{\bex}{\begin{example}}
\nc{\bde}{\begin{definition}}
\nc{\ede}{\end{definition}}
\nc{\epr}{\end{prop}}
\nc{\ethe}{\end{theorem}}
\nc{\ele}{\end{lem}}
\nc{\eco}{\end{corollary}}
\nc{\ere}{\hfill\mbox{$\Diamond$}\end{remark}}
\nc{\eex}{\hfill\mbox{$\Diamond$}\end{example}}
\nc{\epf}{\hfill\mbox{$\Box$}~\\~\\}
\nc{\ba}{\begin{array}}
\nc{\ea}{\end{array}}
\nc{\bea}{\begin{eqnarray}}
\nc{\eea}{\end{eqnarray}}
\nc{\nn}{\nonumber}
\nc{\be}{\begin{enumerate}}
\nc{\ee}{\end{enumerate}}
\nc{\beq}{\begin{equation}}
\nc{\eeq}{\end{equation}}
\nc{\bi}{\begin{itemize}}
\nc{\ei}{\end{itemize}}
\nc{\ra}{\rightarrow}
\nc{\ci}{\circ}
\nc{\lra}{\longrightarrow}
\rnc{\to}{\mapsto}
\nc{\imp}{\Rightarrow}
\rnc{\iff}{\Leftrightarrow}
\rnc{\phi}{\mbox{$\varphi$}}
\rnc{\epsilon}{\varepsilon}
\nc{\al}{\mbox{$\alpha$}}
\nc{\ha}{\mbox{$\alpha$}}
\nc{\hb}{\mbox{$\beta$}}
\nc{\hg}{\mbox{$\gamma$}}
\nc{\hd}{\mbox{$\delta$}}
\nc{\he}{\mbox{$\varepsilon$}}
\nc{\hz}{\mbox{$\zeta$}}
\nc{\hs}{\mbox{$\sigma$}}
\nc{\hk}{\mbox{$\kappa$}}
\nc{\hm}{\mbox{$\mu$}}
\nc{\hn}{\mbox{$\nu$}}
\nc{\la}{\mbox{$\lambda$}}
\nc{\hl}{\mbox{$\lambda$}}
\nc{\hG}{\mbox{$\Gamma$}}
\nc{\hD}{\mbox{$\Delta$}}
\nc{\Th}{\mbox{$\Theta$}}
\nc{\ho}{\mbox{$\omega$}}
\nc{\hO}{\mbox{$\Omega$}}
\nc{\hp}{\mbox{$\pi$}}
\nc{\hP}{\mbox{$\Pi$}}
\nc{\as}{\mbox{$A(S^3\sb s)$}}
\nc{\bs}{\mbox{$A(S^2\sb s)$}}
\nc{\slq}{\mbox{$A(SL\sb q(2))$}}
\nc{\fr}{\mbox{$Fr\llp A(SL(2,\IC))\lrp$}}
\nc{\slc}{\mbox{$A(SL(2,\IC))$}}
\nc{\af}{\mbox{$A(F)$}}
\nc{\suq}{\mbox{$A(SU_q(2))$}}
\nc{\asq}{\mbox{$A(S_q^2)$}}
\nc{\tasq}{\mbox{$\widetilde{A}(S_q^2)$}}
\nc{\Llp}{\mbox{\Large $($}}
\nc{\Lrp}{\mbox{\Large $)$}}
\nc{\llp}{\mbox{\large $($}}
\nc{\lrp}{\mbox{\large $)$}}

\def\ta{\tilde a}
\def\tb{\tilde b}
\def\tc{\tilde c}
\def\td{\tilde d}

\newcommand{\Mat}{{\rm Mat}\,}
\newcommand{\half}{\frac{1}{2}}

\newcommand{\bbr}{{\bf R}}
\newcommand{\bbz}{{\bf Z}}
\newcommand{\Ci}{C_{\infty}}
\newcommand{\Cb}{C_{b}}
\newcommand{\fa}{\forall}

\newcommand{\IC}{{\mathbb C}}
\newcommand{\ID}{{\mathbb D}}
\newcommand{\IF}{{\mathbb F}}
\newcommand{\IH}{{\mathbb H}}
\newcommand{\II}{{\mathbb I}}
\newcommand{\IK}{{\mathbb K}}
\newcommand{\IM}{{\mathbb M}}
\newcommand{\IN}{{\mathbb N}}
\newcommand{\IP}{{\mathbb P}}
\newcommand{\IQ}{{\mathbb Q}}
\newcommand{\IZ}{{\mathbb Z}}
\newcommand{\IR}{{\mathbb R}}
\def\a{\alpha}
\def\b{\beta}
\nc{\ot}{\otimes}
\def\otc{\otimes_{\IC}}
\def\ota{\otimes_ A}
\def\otza{\otimes_{ Z(A)}}
\def\otc{\otimes_{\IC}}
\def\<{\langle}
\def\>{\rangle}
\def\vt{\vartriangleright}
%
%%%%%%%%%%%%%%%%%%%%%%%%%%%%%%%%%%%%%%%%%%%%%%%%%%%%%%%%%%%%%%%%%%%%%%%%%%%%%
%
\title{
\vspace*{1cm} {\large\bf INSTANTONS ON THE QUANTUM $4$-SPHERES $S^4_q$}
\\
\author{
{\sc Ludwik D\c{a}browski}
\vspace*{-1mm}\\
\normalsize Scuola Internazionale Superiore di Studi Avanzati
\vspace*{-1mm}\\
\normalsize Via Beirut 2-4, 34014 Trieste, Italy
\vspace*{-.5mm}\\ dabrow@sissa.it
\vspace*{5mm}\\
{\sc Giovanni Landi}
\vspace*{-1mm}\\
\normalsize Dipartimento di Scienze Matematiche, Universit\`a di Trieste
\vspace*{-1mm}\\
\normalsize Via Valerio 12/b, 34127, Trieste, Italy
\vspace*{-.5mm}\\
landi@dsm.univ.trieste.it
\vspace*{5mm}\\
{\sc Tetsuya Masuda}
\vspace*{-1mm}\\
\normalsize Institute of Mathematics, University of Tsukuba
\vspace*{-1mm}\\
\normalsize
Tsukuba 305, Japan
\vspace*{-.5mm}\\
tetsuya@math.tsukuba.ac.jp}
\date{}
}
\maketitle

\vspace{1cm}
\begin{abstract}
\noindent
We introduce noncommutative algebras $A_q$ of quantum $4$-spheres $S^4_q$,
with $q\in\IR$, defined via a suspension of the quantum group $SU_q(2)$, and a
quantum instanton bundle described by a selfadjoint idempotent $e\in 
\Mat_4(A_q)$,
$e^2=e=e^*$. Contrary to what happens for the classical case or for the
noncommutative instanton constructed in
\cite{cl00}, the first Chern-Connes class $ch_1(e)$ does not vanish thus
signaling a dimension drop.
The second Chern-Connes class $ch_2(e)$ does not vanish as well and the couple
$(ch_1(e), ch_2(e))$ defines a cycle in the $(b,B)$ bicomplex of cyclic
homology.
\end{abstract}

\vfill\eject

\section{Introduction}
The goal of this paper is to provide more interesting examples
of globally nontrivial four dimensional quantum manifolds
and vector bundles (finite projective modules) over them.
Most of the literature concentrated so far on noncommutative tori
or Moyal deformations of $\IR^4$ \cite{rif,rifsch,nesch,furu,KKO}.

In \cite{cl00} an instanton bundle over  a family of noncommutative 4-spheres
$S_{\theta}^4$ which are suspensions of a class of noncommutative 3-spheres,
were introduced.
These spaces $S_{\theta}^4$ fulfill all the axioms
of Riemannian spin geometry as formulated in \cite{co96,co98}.
Moreover, as it happens \cite{cojmp} for the instanton bundle over the ordinary
$4$-sphere $S^4$, the $0$th and the $1$st Chern-Connes classes  of the
instanton
idempotent $e$ vanish, ${\rm ch}_j (e) = 0 \, , ~j = 0, 1$, while
${\rm ch}_2 (e)$ is a nontrivial Hochschild cycle.

In the present paper, we exhibit another family of quantum 4-spheres $S^4_q$,
$q\in \IR$,  defined also via a suspension but now of the quantum 3-sphere
$S^3_q$,  which we take just as the underlying `space'
of the quantum group $SU_q(2)$.
If fact the present spheres were discovered before the ones in
\cite{cl00}.  Though these two families are in a sense
related by analytic continuation of the deformation parameter \cite{dala},
they present some different and worth mentioning properties.
First of all, $S^4_q$ does not seem to obey (all of) the axioms
of a noncommutative manifold as given in \cite{co96,co98},
a feature which is shared with most of the quantum spaces defined
in the framework of the so called $q$-deformations and quantum groups.
Indeed, this fact may be an inspiration for weakening some of those
axioms in order to embrace such a class of spaces.
Furthermore, the quantum spheres $S^4_q$ come equipped with a natural
idempotent $e$ as well which determines a vector bundle (i.e. a finite 
projective module of sections) over it. 
However, contrary to what happens for the 
classical case  or
for the instanton bundle constructed in \cite{cl00}, now the first 
Chern-Connes 
class
$ch_1(e)$ does not vanish. Thus the idempotent $e$ does not provide a 
representation of the
universal instanton algebras as defined in \cite{cojmp,cl00}.
It turns out that the second Chern-Connes class
$ch_2(e)$  does not vanish as well and the couple $(ch_1(e), ch_2(e))$
defines a cycle in the $(b,B)$ bicomplex of cyclic homology \cite{co85,L}.

\section{The algebra of $S^4_q$}

We define the quantum 4-sphere $S^4_q$ as the suspension of the quantum
3-sphere
$S^3_q$  which we take as the underlying `space' of the quantum group
$SU_q(2)$.
Thus, using the definition of the $C^*$-algebra of $S_{\mu}U(2)$ as given in
\cite{SLW-RIMS,SLW-CMP} (with
the convenient replacements
$\mu \mapsto q$, $\alpha \mapsto \alpha^*$ and $\hg \mapsto \hb^*$),
for a parameter 
$q \in [-1,0)\cup (0,1]$ 
we
define $A_q$ as the $C^*$-algebra with unit $\II$ generated by three elements
$\alpha
, \hb $ and $z$ satisfying the relations
\bea\label{s4rel}
&& \hb\alpha = q\alpha\hb , ~~~
\hb^*\alpha = q \alpha\hb^* , ~~~
\hb\hb^* = \hb^*\hb , \\
&& z=z^*, ~~~ z\alpha = \alpha z,~~~ z\hb = \hb z,
~~~
\nonumber\\
&& \alpha^* \alpha + q^2\hb^*\hb +z^2 = \II,
~~~~~~\alpha\alpha^* + \hb\hb^* + z^2 = \nonumber \II ~.
\eea
In particular the `suspension' generator $z$ is central and selfadjoint. It
should be clear from the previous relations that, as it happens for the quantum
group $SU_q(2)$, the generator $\al$ is not normal, $\alpha^* \alpha \not=
\alpha\alpha^*$.

More precisely the algebra $A_q$ can be defined in the following way.
Consider the free (noncommutative) \mbox{$^*$-algebra}
with unity $~\IC [[\alpha , \hb , z, \alpha^* , \hb^* , z^*]]~$
generated by three elements $\alpha , \hb$ and $z$.
A \mbox{$^*$--representation}
$\pi$ of $\IC [[\alpha , \hb , z, \alpha^* , \hb^*, z^*]]$ in terms of
bounded operators on a Hilbert space $H$ is said to be {\em admissible}
if the operators $\pi(\alpha ), \pi(\hb ), \pi(z)$ satisfy the relations
(\ref{s4rel}).
Then for arbitrary
$a\in
\IC [[\alpha , \hb , z, \alpha^* , \hb^*, z^*]]$ we set
$\| a \|$ to be the supremum,
over all admissible representations $\pi$ of
$\IC [[\alpha , \hb , z, \alpha^*, \hb^*, z^*]]$,
of the operator norms $\| \pi (a)\|$.
It can be seen that $\| a \| < \infty$ and that $\| \cdot \|$ is a
$C^*$-semi-norm.
As a consequence the set $\cal J$ of all those
$a$ in $\IC [[\alpha , \hb, z, \alpha^* , \hb^*, z^*]]$ with $\| a \| = 0$
is a two-sided ideal in
$\IC [[\alpha , \hb , z, \alpha^* , \hb^*, z^*]]$. Then one obtains
a $C^*$-norm on the quotient algebra
$\IC [[\alpha , \hb , z, \alpha^* , \hb^*, z^*]]/ \cal J$,
the completion of which is the $C^*$-algebra $A_q$ in question.

In the sequel, to simplify the notation we shall denote the
$\cal J$-equivalence class $a + \cal J$ simply by $a$.
With this convention, one sees that the $^*$-subalgebra generated by $\alpha ,
\hb$ and $z$ is dense in $A_q$.
Moreover, for any triple $\hat{\alpha} , \hat{\hb} ,
\hat{z} $ of bounded operators on $H$
satisfying the relations (\ref{s4rel}) there exists exactly one representation
$\pi : A_q \ra B(H)$ such that $\pi (\alpha ) = \hat{\alpha}$, $\pi (\hb ) =
\hat{\hb}$ and $\pi (z) = \hat{z} $.
It can be also verified that the ideal $\cal J$
actually coincides with the ideal in $\IC [[\alpha , \hb , z, \alpha^* , \hb^*,
z^*]]$ generated by the following elements:
\bea
&& \hb \alpha - q \alpha \hb~, ~~~
\hb^*\alpha - q \alpha\hb^*~, ~~~
\hb^*\hb - \hb\hb^*~,~ \\
&& \alpha^*\hb^* - q\hb^* \alpha^*~, ~~~
\alpha^* \hb- q\hb \alpha^*~,~ \nonumber\\
&& z-z^*,~ z \alpha - \alpha z~, ~~~ z \hb- \hb z~,
\nonumber\\
&& \alpha^* \alpha + q^2\hb^*\hb +z^2 - \II~,
~~~ \alpha\alpha^* + \hb^*\hb +z^2 - \II~. \nonumber
\eea
Using the relations (\ref{s4rel}) a linear basis for $A_q$ can be taken
as $a_{kmn\ell}$, with $k\in\IZ$ and $m, n, \ell$ non negative
integers, of the form
\beq\label{basis}
a_{kmn\ell} = \cases{
\alpha^{*k} \beta^{*m}\beta^n z^{\ell}
~~~~~~\mbox{for}~k=0, 1, 2 ~\dots  \cr
\alpha^{-k} \beta^{*m}\beta^n z^{\ell}
~~~~~~\mbox{for}~k= -1, -2 ~\dots ~. }
\eeq

\noindent
Notice that the quantum sphere $S^4_q$ may be defined also for $q\in\IR$,
$|q| > 1$, but with the
transformation $q \to 1/q$, $\alpha \to \al^*$, $\hb \to q\hb$ and $z \to z$,
we get a sphere which is $C^*$-isomorphic to one for $|q| < 1$. \\
It is clear that the quotient of the $C^*$-algebra $A_q$ by the ideal 
generated by $z$ can be identified with the $C^*$-algebra of the compact
quantum group $SU_q(2)$. However, in this paper we shall not make any use of
additional structures (like coproduct, counit, and
antipode) coming from $SU_q(2)$. 
In \cite{SLW-RIMS} it was shown that for $q \in (-1,0) \cup (0,1)$ the 
spaces $SU_q(2)$ are all homeomorphic in the sense that the 
corresponding $C^*$-algebras are isomorphic. Then, for $q \in (-1,0) \cup 
(0,1)$, all our $C^*$-algebras $A_q$ are isomorphic as well and all 
corresponding spheres are homeomorphic.

\bigskip
For the generic situation when $-1<q<0$ or $0<q<1$ any character $\chi$
of $A_q$ has to satisfy the equations
\bea\label{s4relclq}
&& \chi(\al^*)=\overline{\chi(\al)}~, ~~~ \chi(\hb^*)=\overline{\chi(\hb)}~,
~~~ \chi(z^*)=\chi(z) ~, \\
  &&\chi(\hb) = 0 ~~{\rm and }~~ |\chi(\al)|^2+(\chi(z))^2 = 1 ~. \nn
\eea
To show that the space of all characters is homeomorphic to
the two dimensional sphere $S^2$, we take a generic $\alpha' \in \IC $ and
$z'\in \IR$ such that $|\al^\prime |^2 + z'^{2}=1$. Then, from the general
considerations presented above, there is a 1-dimensional representation
(that is a character) $\chi$ of $A_q$ such that $\chi (\al)=\al'$,
$\chi (\hb )= 0 $ and $\chi (z)= z'$ and this proves the homeomorphism in
question. Hence, for $-1<q<0$ or $0<q<1$ the space of (nonzero) characters
of $A_q$, which can be thought of as the space of `classical points' of
$S^4_q$, is homeomorphic to the classical $S^2$.

For the particular case $q = 1$ the algebra of the sphere $S^4_{q}$ is
commutative. The associated space of characters is
homeomorphic to the $4$-dimensional sphere $S^4$. Indeed any character
$\chi$ of $A_{q=1}$ satisfies the equations
\bea\label{s4relcl}
&& \chi(\al^*)=\overline{\chi(\al)}~, ~~~ \chi(\hb^*)=\overline{\chi(\hb)}~,
~~~ \chi(z^*)=\chi(z) ~, \\
{\rm and } && |\chi(\al)|^2+|\chi(\hb)|^2 +(\chi(z))^2 = 1 ~. \nn
\eea
To show that any element of
$S^4$ arises in this way, similarly to what we did before we take generic
$\alpha' , \hb' \in \IC$ and $z' \in \IR$ such that
$|\al^\prime |^2 + |\hb^\prime |^2 +
z'^{2}=1$. Thus they satisfy relations (\ref{s4relcl}) (or relations
(\ref{s4rel}) for $q=1$) and there is a 1-dimensional representation
$\chi$ of $A_q$ ($q=1$) such that $\chi (\al)=\al '$, $\chi (\hb )=\hb'$ and
$\chi (z)= z'$. This proves the homeomorphism in question and shows that the
algebra $A_q$ for $q=1$ can be identified with the algebra of all continuous
functions on the $4$-dimensional sphere $S^4$.
It is in this sense that $S^4_q$
provides a deformation of the classical $S^4$.

\bigskip
Next, we describe irreducible representations of the algebra $A_q$
(for $-1<q<0$ or $0<q<1$) as
bounded operators on a infinite dimensional Hilbert space $H$ with
an orthonormal basis $\{ \psi_n ~, ~ n = 0, 1, 2, \cdots ~ \}$.
With $\lambda \in \IC~, ~|\lambda| \leq 1$, we get two
families of representations
$\pi_{\lambda, \pm} : A_q \rightarrow B(H)$ given by
\beq\label{reps}
\begin{array}{ll}
\pi_{\lambda, \pm}(z) \, \psi_n = \pi_{\lambda, \pm}(z^*) \, \psi_n = \pm
\sqrt{ 1 - |\lambda|^2 } \, \psi_n ~,
& ~ \\
\pi_{\lambda, \pm}(\alpha) \, \psi_n = \lambda \, \sqrt{ 1 - q^{2(n+1)} } \,
\psi_{n+1} ~,
& ~~~ \pi_{\lambda, \pm}(\alpha^*) \, \psi_n = \bar{\lambda} \,
\sqrt{ 1 - q^{2n} }
\, \psi_{n-1} ~, \\
\pi_{\lambda, \pm}(\beta) \, \psi_n = \lambda \, q^n \, \psi_n ~, &
~~~
\pi_{\lambda, \pm}(\beta^*) \, \psi_n = \bar{\lambda} \, q^n \, \psi_n
~,
\end{array}
\eeq
To be precise, for $\lambda$ such that $|\lambda| = 1$, the two
representations $\pi_{\lambda,+}$ and $\pi_{\lambda,-}$ are identical so
that, in fact,
we  have a family of representations parametrized by points on
a classical sphere $S^2$, similarly to what happens for one dimensional
representations (characters) as described before. \\
As mentioned already, the quotient of the $C^*$-algebra 
$A_q$ by the ideal generated by $z$ is the $C^*$-algebra of the compact
quantum group $SU_q(2)$. Then, with $|\lambda| = 1$, the representations
$\pi_{\lambda,+} = \pi_{\lambda,-} =: \pi_{\lambda}$ yield representations of
$SU_q(2)$ which are unitary equivalent to the ones constructed by Woronowicz
(see for instance (\cite{wdn})).

\section{The instanton and its classes}

Consider now the following element $e$ in the algebra $\Mat_4(A_q) \simeq
\Mat_{4}(\IC ) \otimes A_q$
\beq
\label{proj} e = \half\pmatrix{ 1\! +\! z, ~~0,~~ ~~~\alpha ,~ ~~\hb~ \cr
~~~0,~~
1\! +\! z, ~-\! q\hb^* , ~\alpha^*\cr ~~~\alpha^*, ~-\! q\hb , ~~1\! -\! z,
~~0~~
\cr ~\hb^*,~ ~~~\alpha,~~ ~~~0,~~ 1\! -\! z}\ .
\eeq
Using the relations
(\ref{s4rel}) it can be verified that $e$ is a selfadjoint idempotent
(projection)
\[ e^2 = e = e^*\ .
\] It operates on the right $A_q$-module $A_q^4 = A_q\otimes \IC^4$ and its
range
may be thought of as sections of a vector bundle over $S^4_q$. It is easy to
see that
$e A_q^4$ is a deformation of the classical instanton bundle over $S^4$
in the sense that for $q=1$, the module $e A_q^4$ is the module of sections of
the complex rank two instanton bundle over $S^4$ \cite{at}.

Next, we compute the Chern-Connes Character of the idempotent $e$ given in
(\ref{proj}).  If $\langle \ \rangle$ is
the projection on the commutant of $4 \times 4$ matrices,
up to normalization the component of the (reduced) Chern-Connes Character
are given by
\beq
ch_n(e) = \left\langle \left(e - {1 \over 2} \right) \ot \underbrace{e
\ot \cdots
\ot e}_{2n}
\right\rangle ~, ~~~ n = 0, 1, 2, \dots ~, 
\eeq
and they are elements of
\beq
A_q \ot \underbrace{\bar{A_q} \ot \cdots \ot \bar{A_q}}_{2n} ~,
\eeq
where $\bar{A_q} = A_q / \IC \II$ is the quotient of the algebra $A_q$ by the
scalar multiples of the unit $\II$.\\
The crucial property of the
components ${\rm ch}_n (e)$ is that they define a {\it cycle} in the $(b,B)$
bicomplex of cyclic homology \cite{co85,L}, that is,
\begin{equation}
B \, {\rm ch}_n (e) = b \, {\rm ch}_{n+1} (e) \, . \label{Bb}
\end{equation}
The operator $b$ is defined by
\beq\label{bi}
b (a_0 \ot a_1 \ot \cdots \ot a_m ) =
\sum_{j=0}^{m-1} (-1)^j a_0 \ot \cdots \ot a_j a_{j+1} \ot \cdots \ot a_m
\, + (-1)^{m} a_m a_0 \ot a_1 \ot \cdots \ot a_{m-1} \,
\eeq
while the operator $B$ is written as
\beq\label{Bi}
B = A B_0 \, ,
\eeq
where
\bea\label{Bi0}
&& B_0 (a_0 \ot a_1 \ot \cdots \ot a_m ) = \II \ot a_0 \ot a_1 \ot \cdots
\ot a_m
\\
&& A (a_0 \ot a_1 \ot \cdots \ot a_m ) = { 1 \over m} \sum_{j=0}^m 
(-1)^{mj} a_j \ot
a_{j+1} \ot \cdots \ot a_{j-1} \, , \label{Ai}
\eea
with the obvious cyclic identification $m+1 = 0$.
To be precise, in formul{\ae} (\ref{bi}), (\ref{Bi0}) and (\ref{Ai}),
all elements in the tensor products but the first one should be taken modulo
complex multiples of the unit $\II$, that is one has to project onto
$\bar{A_q} = A_q / \IC \II$.

\bigskip
For the $0$th component of the Chern-Connes Character of the idempotent 
(\ref{proj}) on the spheres $S^4_q$ we find,
\beq
ch_0(e) = \left\langle \left(e - {1 \over 2} \right) \right\rangle = 0 ~.
\eeq
This could be interpreted as saying that the idempotent and the corresponding
module (the `vector bundle') has complex rank equal to $2$.

\bigskip
Next for the $1$st component we have,
\bea
ch_1(e) &=& \left\langle \left(e - {1 \over 2} \right) \ot e \ot e
\right\rangle
\\ &=& {1 \over 8} (1-q^2) \Big\{
z \ot (\beta \ot \beta^* - \beta^* \ot \beta ) \nn \\
&& ~~~~~~~~~~~~~~~~
+ \beta^* \ot (z \ot \beta - \beta \ot z) + \beta \ot (\beta^* \ot z -
z \ot \beta^*) \Big\}
  ~. \nn
\eea
It is straightforward to check that
\beq
b ch_1(e) = 0 = B ch_0(e)
\eeq

\bigskip
Finally, the $2$nd component
\beq
ch_2(e) = \left\langle \left(e - {1 \over 2} \right) \ot e \ot e\ot e \ot
e \right\rangle
\eeq
can be written as a sum of five terms
\beq
ch_2(e) = {1\over 32} \, \Big( z \, \ot \, c_z + \a \, \ot \, c_{\a} + \a^*
\, \ot \,
c_{\a^*} +
\b \, \ot \, c_{\b} + \b^* \, \ot \, c_{\b^*} \Big) \, , \label{hocy}
\eeq
with
\begin{eqnarray}\label{gamz}
c_z &=& (1-q^4) \, (\b \, \ot \, \b^* \, \ot \, \b \, \ot \, \b^*
-  \b^* \, \ot \, \b \, \ot \, \b^* \, \ot \, \b) \\
&+& (1-q^2) \, \Big\{
z \, \ot \, z \, \ot \, (\b \, \ot \, \b^* -  \b^* \, \ot \, \b)
+ (\b \, \ot \, z \, \ot \, z \, \ot \, \b^* -  \b^*  \, \ot \, z \, \ot \,
z \, \ot
\, \b) \nn \\
&\, & \, \, \, \, \, \, \, \, \, \, \, \, \, \, \, \,\, \, \,\, \, \,\, \, \,
+  \, (\b \, \ot \, \b^* - \b^* \, \ot \, \b)  \, \ot \,
z \, \ot \, z  + z \, \ot \, (\b \, \ot \, \b^* - \b^* \, \ot \, \b) \, \ot
\, z \nn
\\ &\, & \, \, \, \, \, \, \, \, \, \, \, \,
- z \, \ot \, (\b \, \ot \, z \, \ot \, \b^* -  \b^* \, \ot \, z \, \ot \, \b)
- (\b \, \ot \, z
\, \ot \, \b^* -  \b^* \, \ot \, z \, \ot \, \b) \, \ot \, z  \, \Big\}
\nn \\
&+& (\a \, \ot \, \a^* - q^2 \, \a^* \, \ot \, \a)
\, \ot \, (\b \, \ot \, \b^* -  \b^* \, \ot \, \b)
\nonumber \\ &\, & \, \, \, \, \, \, \, \, \, \,
+ \, ( \b \, \ot \, \b^* -  \b^* \, \ot \, \b)
\, \ot \, ( \a \, \ot \, \a^* - q^2 \, \a^* \, \ot \, \a)
\nonumber \\
&+& \, (\b \, \ot \, \a  - q \, \a \, \ot \, \b)
\, \ot \, (\a^* \, \ot \, \b^* - q \, \b^* \, \ot \, \a^*)
\nonumber \\
&\, & \, \, \, \, \, \, \, \, \, \,
+ \, (\a^* \, \ot \, \b^* - q \, \b^* \, \ot \, \a^*)
\, \ot \, (\b \, \ot \, \a  - q \, \a \, \ot \, \b)
\nonumber \\
&+& \, ( \a^* \, \ot \, \b - q \,  \b \, \ot \, \a^*)
\, \ot \, (q \,  \a \, \ot \, \b^* \, -  \b^* \, \ot \, \a)
\nonumber \\
&\, & \, \, \, \, \, \, \, \, \, \,
+ \, (q \,  \a \, \ot \, \b^* \, -   \b^* \, \ot \, \a)
\, \ot \, ( \a^* \, \ot \, \b - q \,  \b \, \ot \, \a^*) \, ; \nonumber
\end{eqnarray}

\begin{eqnarray}\label{gamalp}
c_\a &=&
(z \, \ot \, \a^* - \a^* \, \ot \, z)
\, \ot \, (\b^* \, \ot \, \b -  \b \, \ot \, \b^*)
\\ &\, & \, \, \, \, \, \, \, \, \, \,
+ \, q^2 \, (\b^* \, \ot \, \b -  \b \, \ot \, \b^*)
\, \ot \, (z \, \ot \, \a^* - \a^* \, \ot \, z)
\nonumber \\
&+& \, q \, (z \, \ot \, \b  - \b \, \ot \, z)
\, \ot \, (\a^* \, \ot \, \b^* - q \, \b^* \, \ot \, \a^*)
\nonumber \\
&\, & \, \, \, \, \, \, \, \, \, \,
+ \, (\a^* \, \ot \, \b^* - q \, \b^* \, \ot \, \a^*)
\, \ot \, (z \, \ot \, \b  - \b \, \ot \, z)
\nonumber \\
&+& \, q \, (\b^* \, \ot \, z \, - z \, \ot \, \b^*)
\, \ot \, ( \a^* \, \ot \, \b - q \,  \b \, \ot \, \a^*)
\nonumber \\
&\, & \, \, \, \, \, \, \, \, \, \,
+ \, ( \a^* \, \ot \, \b - q \, \b \, \ot \, \a^*)
\, \ot \, (\b^* \, \ot \, z \, - z \, \ot \, \b^*) \, ; \nonumber
\end{eqnarray}
\begin{eqnarray}\label{gamalpbar}
c_{\a^*} &=&
q^2 \, (z \, \ot \, \a - \a \, \ot \, z)
\, \ot \, (\b \, \ot \, \b^* -  \b^* \, \ot \, \b)
\\ &\, & \, \, \, \, \, \, \, \, \, \,
+ \, (\b \, \ot \, \b^* -  \b^* \, \ot \, \b)
\, \ot \, (z \, \ot \, \a - \a \, \ot \, z)
\nonumber \\
&+& \, (\b^* \, \ot \, z - z \, \ot \, \b^* )
\, \ot \, (\b \, \ot \, \a - q \, \a \, \ot \, \b)
\nonumber \\
&\, & \, \, \, \, \, \, \, \, \, \,
+ \, q \, (\b \, \ot \, \a - q \, \a \, \ot \, \b)
\, \ot \, (\b^* \, \ot \, z - z \, \ot \, \b^* )
\nonumber \\
&+& \, (z \, \ot \, \b  \, - \b \, \ot \, z)
\, \ot \, (\b^* \, \ot \, \a - q \,  \a \, \ot \, \b^*)
\nonumber \\
&\, & \, \, \, \, \, \, \, \, \, \,
+ \, q \, (\b^* \, \ot \, \a - q \,  \a \, \ot \, \b^*)
\, \ot \, (z \, \ot \, \b  \, - \b \, \ot \, z) \, ; \nonumber
\end{eqnarray}

\begin{eqnarray}\label{gambet}
c_\b &=&
(1-q^4) \, \Big[ (\b^* \, \ot \, z - z \, \ot \, \b^*)  \, \ot \, \b \, \ot
\, \b^*
+  \b^* \, \ot \, \b \, \ot (\b^* \, \ot \, z - z \, \ot \, \b^*) \Big] \\
&+& (1-q^2) \, \Big\{
\b^* \, \ot \, z \, \ot \, z \, \ot \, z
- z \, \ot \, \b^* \, \ot \, z \, \ot \, z  \nn \\
&\, & \, \, \, \, \, \, \, \, \, \, \, \, \, \, \, \, \, \, \, \,
\, \, \, \, \, \, \, \, \, \, \, \, \, \, \, \, \, \, \, \,
\, \, \, \, \, \, \, \, \, \, \, \, \, \, \, \, \, \, \, \,
+ z \, \ot \, z \, \ot \, \b^* \, \ot \, z
- z \, \ot \, z \, \ot \, z\ot \, \, \b^* \,
  \Big\} \nn \\
&+& (\b^* \, \ot \, z - z \, \ot \, \b^*)
\, \ot \, (\a \, \ot \, \a^* - q^2 \, \a^* \, \ot \, \a)
\nn \\ &\, & \, \, \, \, \, \, \, \, \, \,
+ \, (\a \, \ot \, \a^* - q^2 \, \a^* \, \ot \, \a)
\, \ot \, (\b^* \, \ot \, z - z \, \ot \, \b^*)
\nonumber \\
&+& \, (\a \, \ot \, z - z \, \ot \, \a)
\, \ot \, (\a^* \, \ot \, \b^* - q \, \b^* \, \ot \, \a^*)
\nonumber \\
&\, & \, \, \, \, \, \, \, \, \, \,
+ \, q \, (\a^* \, \ot \, \b^* - q \, \b^* \, \ot \, \a^*)
\, \ot \, (\a \, \ot \, z - z \, \ot \, \a)
\nonumber \\
&+& \, (\b^* \, \ot \, \a - q \, \a \, \ot \, \b^*)
\, \ot \, (\a^* \, \ot \, z \, - z \, \ot \, \a^*)
\nonumber \\
&\, & \, \, \, \, \, \, \, \, \, \,
+ \, q \, (\a^* \, \ot \, z \, - z \, \ot \, \a^*) \, \ot \, (\b^* \, \ot
\, \a - q
\, \a \,
\ot \, \b^*) \, ; \nonumber
\end{eqnarray}

\begin{eqnarray}\label{gambetbar}
c_{\b^*} &=&
(1-q^4) \, \Big[ (z \, \ot \, \b - \b \, \ot \, z)  \, \ot \, \b^* \, \ot
\, \b
+  \b \, \ot \, \b^* \, \ot (z \, \ot \, \b - \b \, \ot \, z ) \Big] \\
&+& (1-q^2) \,
\Big\{ - \b \, \ot \, z \, \ot \, z \, \ot \, z
+ z \, \ot \, \b \, \ot \, z \, \ot \, z
\nonumber \\
&\, & \, \, \, \, \, \, \, \, \, \, \, \, \, \, \, \, \, \, \, \, \, \, \,
\, \, \,
\, \, \, \,\, \, \, \, \, \, \, \, \, \,
- z \, \ot \, z \, \ot \, \b \, \ot \, z
+ z \, \ot \, z \, \ot \, z\ot \, \, \b \,
  \Big\} \nn \\
&+& (z \, \ot \, \b - \b \, \ot \, z)
\, \ot \, (\a \, \ot \, \a^* - q^2 \, \a^* \, \ot \, \a)
\nn \\ &\, & \, \, \, \, \, \, \, \, \, \,
+ \, (\a \, \ot \, \a^* - q^2 \, \a^* \, \ot \, \a)
\, \ot \, (z \, \ot \, \b - \b \, \ot \, z)
\nonumber \\
&+& \, q \, (z \, \ot \, \a^* - \a^* \, \ot \, z)
\, \ot \, (\b \, \ot \, \a - q \, \a \, \ot \, \b)
\nonumber \\
&\, & \, \, \, \, \, \, \, \, \, \,
+ \, (\b \, \ot \, \a - q \, \a \, \ot \, \b)
\, \ot \, (z \, \ot \, \a^* - \a^* \, \ot \, z)
\nonumber \\
&+& q \, (\a^* \, \ot \, \b - q \, \b \, \ot \, \a^*)
\, \ot \, (z \, \ot \, \a - \a \, \ot \, z) \,
\nonumber \\
&\, & \, \, \, \, \, \, \, \, \, \,
+ \, \, (z \, \ot \, \a - \a \, \ot \, z)
\, \ot \, (\a^* \, \ot \, \b - q \, \b \, \ot \, \a^*) \, . \nonumber
\end{eqnarray}

\bigskip
\noindent
By using the relations (\ref{s4rel}) for our algebra, and remembering
that we need to project on $\bar{A_q}$ in all terms of the tensor product
but the
first one, a long (one needs to compute 750 terms) but straightforward
computation gives
\bea\label{bch2}
&& b ch_2(e) = {1 \over 16} \, (1-q^2) \, \Big\{
\, \II \, \ot \, z \, \ot \, ( \b \, \ot \, \b^*
- \b^* \, \ot \, \b )
\\ && \, \, \, \, \, \, \, \, \, \, \, \, \, \, \, \, \, \, \, \, \, \, \,
\, \,
\,\, \, \, \, \, \, \, \, \, \, \, \, \, + \, \II \, \ot \, \b \, \ot \, (\b^*
\ot \, z \, - z \ot \b^*)
  + \, \II \,
\ot \, \b^* \ot \, (z \, \ot \, \b - \b \ot \, z ) \, \Big\} \nn
\eea
and this is exactly equal to $B ch_1(e)$.

\section{Final remarks}

There are several directions in which one can proceed and we just mention some
of them.

\noindent
It would be clearly very interesting to study differential calculi on
our quantum $4$-sphere and develop Yang-Mills theory.

\noindent
Another natural question is to which extent the sphere
$S_{q}^4$ could be endowed with a structure of a metric noncommutative
manifold which fulfills (some of) the related axioms \cite{co96,co98}.
In particular one should construct an appropriate Dirac operator. This will
probably be possible along the lines of \cite{cl00} where it was suggested
that the true Dirac operator $D$ for the quantum $SU_q(2)$ (and also for 
the quantum Podle\'{s} $2$-sphere $S_{q}^2$ \cite{po}) should satisfy 
an equation of the form
\begin{equation}
{ q^{2D} - q^{-2D} \over q^2 - q^{-2}} = Q \, . \label{truedir}
\end{equation}
where $Q$ is some $q$-analogue of the Dirac operator like the ones 
found in \cite{biku,ma00}. \\
Once found the operator $D$, one would easily `suspend' it to the 
4-sphere $S_{q}^4$.

\noindent
Finally, we mention that it will be interesting to study if there is
any relation with the sheaf-theoretic construction of a $q$-deformed 
instanton in \cite{pfl}.

\bigskip\bigskip
\subsection*{Acknowledgements}
We are grateful to Alain Connes for several enlightening conversations.
This work has been partially supported by the Regione   
Friuli-Venezia-Giulia via the Research Project `Noncommutative geometry: 
algebraic, analytical and probabilistic aspects and applications to 
mathematical physics'.

%\vfill\eject

\end{document}